\documentclass[12pt]{article}

\usepackage{fullpage,amsfonts,amsmath,amsthm,amssymb,mathrsfs,graphicx}

\usepackage{verbatim,url}

\usepackage[top=2.54cm, bottom=2.54cm, left=2.54 cm, right=2.54 cm]{geometry}
\newcommand{\lab}[1]{\label{#1}}                

 \usepackage[usenames,dvipsnames]{color}

\newcommand{\remove}[1]{}

\newcommand{\be}{\begin{equation}}
\newcommand{\bel}[1]{\begin{equation}\lab{#1}\ }
\newcommand{\ee}{\end{equation}}
\newcommand{\bea}{\begin{align}}
\newcommand{\eea}{\end{align}}
\newcommand{\bean}{\begin{align*}}
\newcommand{\eean}{\end{align*}}

\newtheorem{thm}{Theorem}

\newtheorem{con}[thm]{Conjecture}

\renewcommand{\ge}{\geqslant}

\def\cjref#1{Conjecture~$\ref{#1}$}

\def\fref#1{Figure~$\ref{#1}$}

\date{}

\title{
Counterexamples on  matchings in hypergraphs and full rainbow matchings in graphs
}

\author{Pu Gao\thanks{Research supported by the ARC grant DE170100716.} \\
{\small School of Mathematical Sciences}\\[-0.8ex]
{\small Monash University}\\
{\small \texttt{jane.gao@monash.edu}} \and 
Reshma Ramadurai\\
{\small School of Mathematics \& Statistics}\\[-0.8ex]
{\small Victoria University of Wellington}\\
{\small \texttt{Reshma.ramadurai@vuw.ac.nz}} \and 
Ian M.\ Wanless\thanks{Research supported by the ARC grant DP150100506.}\\
{\small School of Mathematical Sciences}\\[-0.8ex]
{\small Monash University}\\
{\small \texttt{ian.wanless@monash.edu}} \and 
Nick Wormald \thanks{Research supported by the Australian Laureate Fellowships grant FL120100125.}\\
{\small School of Mathematical Sciences}\\[-0.8ex]
{\small Monash University}\\
{\small \texttt{nick.wormald@monash.edu}}
}

\begin{document}
\maketitle
%

A {graph}  $G$ whose edges are  coloured (not necessarily properly) contains a {\em full rainbow matching} if there is a matching $M$ that contains exactly one edge of each colour. We refute several  conjectures on matchings including the following, {which is equivalent to a weakening of  Conjecture~\ref{cj:falseconj} below, which is in turn due to Aharoni and Berger~\cite[Conj.~2.5]{ahbe}.} 
\begin{con}\label{cj:falseconjbip}
  Let $G$ be a bipartite {graph}  with maximum degree $\Delta(G)$,
  whose edges are  
  coloured {(not necessarily properly)}. If every colour appears on at least $\Delta(G)+1$ edges, then
  $G$ has a full rainbow matching. 
\end{con}
To be precise, we disprove Conjectures 2.5 and 2.9
in \cite{ahbe},  as well as Conjectures 5.3, 5.4, 6.1 and 6.2 made by
  Aharoni, Berger, Chudnovsky,   Howard and  Seymour~\cite{ABCHS}.  

To discuss these, we convert a {graph}  $G$  to a 3-uniform hypergraph $H$ as follows.
The vertices of $H$ are $V(G)\cup V_1$ where $V_1$ is the set of colours
used on edges of $G$. For each edge $\{u,v\}$ of $G$ with colour $c\in V_1$
there is a hyperedge $\{u,v,c\}$ in $H$.
Now a full rainbow matching in $G$ corresponds to a matching of
$H$ that covers all of the vertices in $V_1$. We call this a {\em $V_1$-matching}.  If $G$ happens to be
bipartite with bipartition $V_2\cup V_3$,
then $H$ will be {\em tripartite}, because its vertices can
be partitioned as $V_1\cup V_2\cup V_3$ such that every hyperedge
includes one vertex from each of these three sets. 

\def\GG{\mathcal{G}}
\def\HH{\mathcal{H}}

\begin{figure}
  \[
  \includegraphics[scale=0.5]{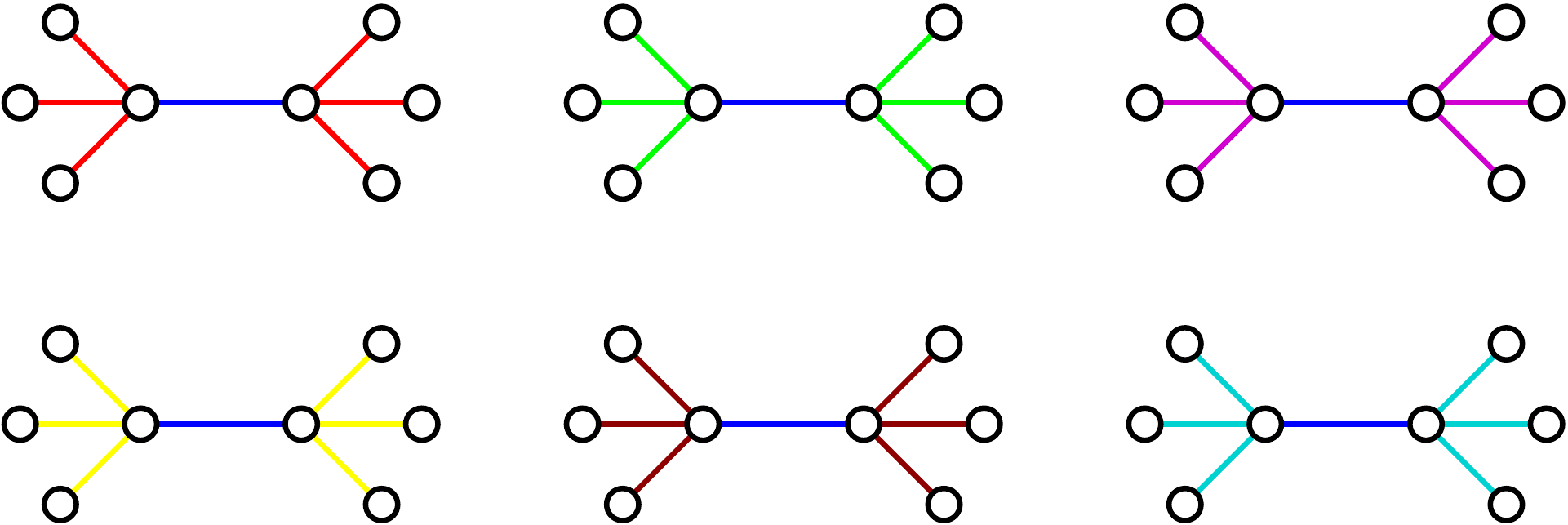}  
  \]
  \caption{\label{f:dblstar}The graph $\GG_6$}
\end{figure}

Let $m$ be a positive even integer.
We now construct a bipartite graph $\GG_m$  that provides a counterexample to Conjecture~\ref{cj:falseconjbip}. 
 There are $m$ components in $\GG_m$,
each isomorphic to a double star which has two adjacent central
vertices each of which has $m/2$ leaves attached to it. The edge
between the central vertices in each double star is coloured blue. In
each component, the edges  incident with  leaves all have one
colour (not blue), which is specific to that component. Hence there are
$m+1$ colours overall, and each colour appears on $m$ edges.
\fref{f:dblstar} shows $\GG_6$.
There is
no full rainbow matching in $\GG_m$ because such a matching must include a blue
edge from some double star $S$. However, the colour of the other edges
in $S$ then cannot be represented in the matching.

Let $\HH_m$ be the tripartite hypergraph corresponding to $\GG_m$.
Let $V_1$ be the vertices of
$\HH$ corresponding to the colours, and $V_2,V_3$ the sets of vertices
corresponding to a bipartition of $\GG$.  Then every vertex in $V_1$ has
degree $m$. The vertices in $V_2\cup V_3$ all have degree either $1$
or $m/2+1$.  For all even $m\ge 4$, this disproves the following conjecture of Aharoni and Berger 
\cite[Conj.~2.5]{ahbe} (repeated as \cite[Conj.~5.3]{ABCHS}).

\begin{con}\label{cj:falseconj}
  Let $H$ be a hypergraph with a vertex tripartition $V(H)=V_1\cup V_2\cup V_3$
  such that every hyperedge includes exactly one vertex from $V_i$ for $i=1,2,3$.
  If $\delta(V_1)>\Delta(V_2\cup V_3)$ then $H$ has a $|V_1|$-matching.
\end{con}

{Another way of saying the above is that  Conjecture~\ref{cj:falseconjbip} is equivalent to   Conjecture~\ref{cj:falseconj} with an appropriate restriction on $H$.}

\begin{figure}
  \[
  \includegraphics[scale=0.4]{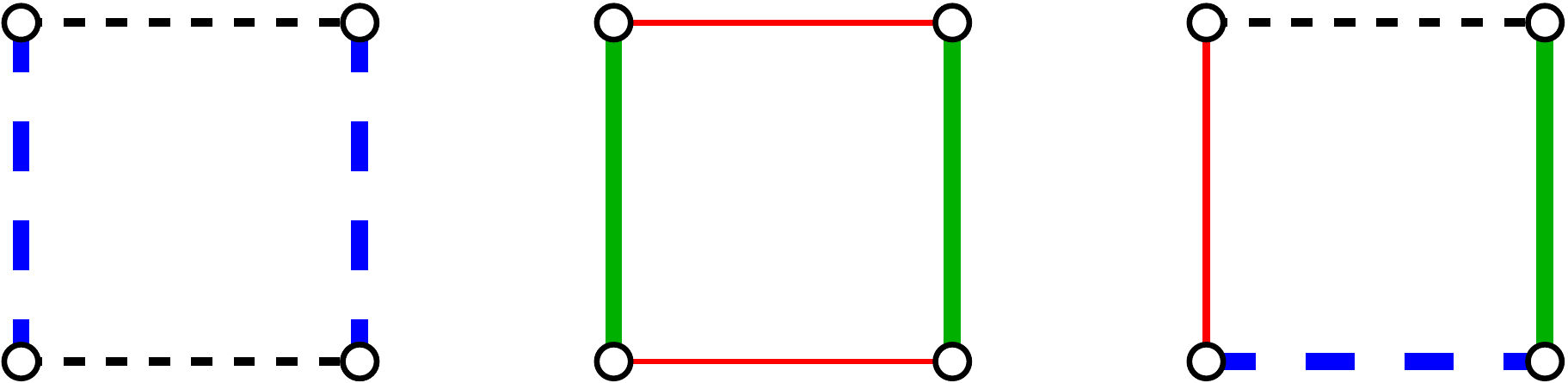}  
  \]
  \caption{\label{f:4cycleCE}A $2$-regular graph with no rainbow matching.}
\end{figure}

 We offer another (single) counterexample to \cjref{cj:falseconj}
  based on the graph  in \fref{f:4cycleCE}, which has no rainbow matching.
The corresponding tripartite hypergraph has 
$\delta(V_1)=3>2=\Delta(V_2\cup V_3)$. The line graph of the graph
in \fref{f:4cycleCE} was  described  in \cite{Alon92} and its complement
figured  in \cite{Jin92}.   In both cases the focus of the investigation
was slightly different from ours, so the generalisations that were offered
are not relevant for us.

  Interestingly, Aharoni and Berger
\cite[Thm~2.6]{ahbe} showed that in any tripartite hypergraph, if
$\delta(V_1)\ge2\Delta(V_2\cup V_3)$ then there must be a
$|V_1|$-matching. Our hypergraph $\HH_m$ shows that their theorem is close
to tight, since the (minimum) degree $\delta(V_1)$ of a vertex in
$V_1$ is nearly double the maximum degree $\Delta(V_2\cup V_3)$ of the
vertices outside $V_1$. 

Conjecture 2.9 of \cite{ahbe} generalises \cjref{cj:falseconj}, so it too is false.
Similarly, \cite[Conj.\,6.1]{ABCHS} asserts that
if $\delta(V_1)\ge2+\Delta(V_2\cup V_3)$ then there must be a
$|V_1|$-matching, so $\HH_m$ is a counterexample whenever $m\ge6$.

Finally, we consider Conjectures 5.4 and 6.2 from
\cite{ABCHS}. These {assert that a full rainbow matching exists
in an edge-coloured graph $G$ (not necessarily bipartite) provided that the number of edges of each colour is at least $\Delta(G)+2$}. Again, $\HH_m$ provides a counterexample. Indeed, it
shows that the $2$ cannot be replaced by any constant.

  \let\oldthebibliography=\thebibliography
  \let\endoldthebibliography=\endthebibliography
  \renewenvironment{thebibliography}[1]{%
    \begin{oldthebibliography}{#1}%
      \setlength{\parskip}{0.4ex plus 0.1ex minus 0.1ex}%
      \setlength{\itemsep}{0.4ex plus 0.1ex minus 0.1ex}%
  }%
  {%
    \end{oldthebibliography}%
  }

\end{document}